\documentclass[12pt,reqno]{amsart}
\usepackage{amsthm, amsfonts,  amsmath, amssymb}
\usepackage{url}
 2

\newtheorem{thm}{Theorem}[section]
\newtheorem{lem}[thm]{Lemma}
\newtheorem{prop}[thm]{Proposition}

\newtheorem{conj}[thm]{Conjecture}

\newtheorem{cor}[thm]{Corollary}
\numberwithin{equation}{section}

\newcommand{\N}{\mathbb{N}}

\newcommand{\BigO}[1]{\ensuremath{\operatorname{O}\bigl(#1\bigr)}}

\title[On sign changes for almost prime coefficients]{On sign changes for almost prime coefficients of half-integral weight modular forms}
\author[S. Krishnamoorthy]{Srilakshmi Krishnamoorthy}
\email{srilakshmi@iitm.ac.in}
\address{Department of Mathematics, Indian Institute of Technology Madras\\
Chennai 600036, INDIA}

\author[M. R. Murty]{M. Ram Murty}
\email{murty@mast.queensu.ca}
\address{Department of Mathematics \\
Queen's University \\
Kingston\\
Ontario \\
K7L 3N6, CANADA.}
            
\begin{document}
\begin{abstract}
For a half-integral weight modular form $f = \sum_{n=1}^{\infty} a_f(n)n^{\frac{k-1}{2}} q^n$ of weight $k=\ell+\frac{1}{2}$ on $\Gamma_0(4)$
such that $a_f(n)$($n \in \N$) are real, we prove for a fixed suitable natural
number $r$ that $a_f(n)$ changes sign infinitely often as $n$ varies over numbers having at most
 $r$ prime factors,
assuming the analog of the Ramanujan conjecture
for Fourier coefficients of half-integral weight forms.
\end{abstract}
\maketitle

\footnotetext[1]{Research of the first author was supported by an DST-INSPIRE Grant. }

\footnotetext[2]{Research of the second author was supported by an NSERC Discovery Grant.}

\section{Introduction and Statement of Results}

Let $f = \sum_{n=1}^{\infty} a_f(n) n^{\frac{k-1}{2}} q^n$ be a half-integral weight modular form in the Kohnen's $+$ subspace of weight $k =\ell+ \frac{1}{2}$ on $\Gamma_0(4)$  with $\ell \geq 2$.
Throughout the article, we shall assume that $a_f(n)$'s are real. 
In~\cite{Koh10}, Kohnen proved that for any half-integral weight modular form $f$, not necessarily an eigenform, for a square-free natural number $t$, the sequence $a_f(tn^2)(n\in \N)$ changes sign infinitely often, provided there exists $n_0$ such that
$a_f(tn^2_0) \neq 0.$ In \cite{HKKL12}, Hulse, Kiral, Kuan and Lim proved that if $f$ is an eigenform, then the sequence $a_f(t)$, where $t$ runs over all square-free integers, changes sign infinitely often. 
In ~\cite{MM13}, Meher and Murty obtained some quantitative results on the number of sign changes in the sequence $a_f(n) (n \in \N).$

Motivated by the Sato-Tate equidistribution theorem for integral weight Hecke eigenforms, it is natural to ask whether there is an equidistribution theorem in the half-integral weight setting also.
To understand this question, the first step is to understand if the sequence $a_f(p)$, where $p$ runs over prime numbers, changes sign infinitely often or not. Though we are not completely successful in answering
this question, we try to answer this question to some extent in this paper.
To this end, we consider the set ${\mathcal P}_r$ of numbers having at most $r$
prime factors. 
Elements of ${\mathcal P}_r$ are called almost primes.
 We sometimes say $n$ is $P_r$
or write $n=P_r$ if $n \in {\mathcal P}_r$.
In this context, this paper will focus on the sequence $a_f(n)$ where $n$ varies over ${\mathcal P}_r$ for a suitable fixed $r$ (to be specified later).
The purpose of this paper is to prove the following theorem:

\begin{thm}\label{mainthm1}
Let $\ell \geq 2$ be an integer and let $f$ be an eigenform in the Kohnen + 
subspace of weight $k= \ell +1/2$ on $\Gamma_0(4)$.  Write
 $$f(z) = \sum^{\infty}_{n=1} a_f(n) n^{\frac{k-1}{2}} q^n . \qquad q=e^{2\pi i z},$$
and assume
\begin{enumerate}
\item  $a_f(n)$ are real for all $n \geq 1$,
\item the Ramanujan conjecture holds for $f$; that is, for all $\epsilon >0$,
and all $n \geq 1$, we have $a_f(n) \ll_{\epsilon} n^{\epsilon}$.
\end{enumerate}
Then, there is an $r \geq 1$ such that the sequence of numbers $a_f(n)$
changes sign infinitely often as $n$ ranges over numbers in  ${\mathcal P}_r$.  
More precisely, the number
of sign changes with $n=P_r$ and $n\leq x$ is $\gg \log x $.
\end{thm}

In~\cite{Mur83}, the oscillations of Fourier coefficients of normalized Hecke eigenforms of integral weight were studied. However, we can't apply those techniques to half-integral weight forms, 
since there is no multiplicative theory of eigenforms in this setting. In spite of this, we prove Theorem~\ref{mainthm1} by using sieve theoretic techniques.
Of course, we would like to prove the theorem with $r=1$, but this
goal seems out of reach with present knowledge and our sieve technique.

Here is our strategy.  As outlined in \cite{MM13}, we need to have
three ingredients to deduce sign change results for any given sequence
of numbers $a(n)$.  First, we need an estimate of the form
$a(n) = O(n^\alpha)$.  Second, we need an estimate 
$$\sum_{n\leq x} a(n) = O(x^\beta). $$
Third, we need an asymptotic formula 
$$\sum_{n\leq x} a(n)^2 = cx + O(x^\gamma)$$
with $\alpha, \beta, \gamma$ non-negative constants and $c>0$.
Then, if $\alpha+\beta < \gamma < 1$, the sequence has infinitely many sign
changes. 
If $\max( \alpha + \beta, \gamma) < b < 1$, then the last condition can be relaxed to
$$\sum_{x \leq n\leq x +y} a(n)^2 \gg y \ \mathrm{for} \ \mathrm{any} \ y > x^{b} . $$
 In our situation, the assumption of the Ramanujan
conjecture (see section 2 below) gives the first estimate
for our sequence $a_f(n)$ with $n$ an almost prime.  
This can be relaxed a bit.  We make some comments in this context at the end
of the paper.
For the second condition, we modify a result of Duke and Iwaniec
\cite{DI90} who showed the required estimate when the argument is
restricted over primes.  We derive the corresponding result for
almost primes.  We prove:

\begin{prop}\label{lem1}
If $a_f(n)$'s are as in Theorem \ref{mainthm1}, then for any natural number $r$ and for sufficiently
large $x,$ we have for any $\epsilon > 0,$
\begin{equation}\tag{5}{\label{eqn5}}
\sum_{1 \leq n \leq x,\,  n = P_r} a_f(n) = \BigO{x^{\frac{155}{156} +\epsilon}}.
\end{equation}
\end{prop}

For the third condition, we apply a lower bound
sieve technique following a method of Hoffstein and Luo \cite{HL97}.  We show:

\begin{prop}\label{mainthm2}
Let $f$ be as in Theorem~\ref{mainthm1}. 
Then there exists a natural number $r$ 
and $\delta >0$ such that for any $Y >3$, we have 
$$\sum_{n = P_r, Y^\delta < n < Y} a_f^2(n) \gg  \frac{Y}{\log Y}.$$ 
\end{prop}

With these results in place, we derive our main theorem following the axiomatic outline given
above.

\section{Notations and Preliminaries}
For the sake of completeness, 
we review some rudimentary facts about half-integral weight modular forms as well
as highlight why one expects Ramanujan's conjecture to be true. Let $f = \sum_{n \geq 1} a_f(n)n^{\frac{k-1}{2}}q^n$ be a cusp form of weight $k = \ell + \frac{1}{2}$ on $\Gamma_0(4)$. Consider the Dirichlet series 
$$L(f,s):=\sum_{n \geq 1} \frac{a_f(n)}{n^s}$$  which converges for $\Re(s)> \frac{3}{2}$ and represents a holomorphic function in this domain. Then by means of the usual Mellin formula and using some standard arguments one proves that
the series admits an analytic continuation and functional equation for all complex values of $s$.
The reader may consult ~\cite[\S 5]{Shi73} or ~\cite[p. 429]{Koh92}.

From the results of Waldspurger ~\cite{Wa81} (see also \cite{kz}), we know that for all square-free $m$,  $a^2_f(m) =O( L( \frac{1}{2}, g, \chi_m))$, where $\chi_m$ is the quadratic character
$( \frac{ m (-1)^{ \frac{l-1}{2}   }  } {  \cdotp}  )$ and $g = \sum_{n \geq 1} c_g(n)n^{(2l-1)/2}q^n$ is the classical modular form of weight $2\ell$ and level $2$ corresponding to the half-integral modular form $f$ via Shimura's correspondence ~\cite{Shi73}. 
In~\cite{Iwa87}, Iwaniec showed that $a_f(p)$'s are bounded by the factor $p^{\frac{3}{14}}.$ 
The exponent was later improved by Blomer and Harcos to $\frac{3}{16}+\epsilon$
in \cite{blomer}.  
One expects a stronger
bound.  Indeed, if we assume the analog of the Lindel\"of hypothesis for $L(1/2, g, \chi_m)$ in the conductor aspect
(that is, the $m$ -aspect), then it is reasonable to expect the following:
\begin{conj}(Ramanujan conjecture)\label{conj} 
Let  $f(z) = \sum^{\infty}_{n=1} a_f(n) n^{\frac{k-1}{2}} q^n$ be a half-integral weight modular form of weight 
$k=\ell+\frac{1}{2}$ on $\Gamma_0(4N)$,
where $k \in \N,$ $k \geq 2$ and $q = e^{2 \pi i z}.$ Then, for any $\epsilon$ $>$ $0$,
\begin{equation}\tag{1}\label{eqn1}
a_f(n) =  \BigO{n^{\epsilon}}.  
\end{equation}
\end{conj}

We also need the following result of  Duke-Iwaniec (see Section 8 of ~\cite{DI90}).
\begin{prop}\label{lem11}
 If   $a_f(p)$'s  are as in Theorem \ref{mainthm1}, then for sufficiently
 large $x,$ and for any $\epsilon$ $>$ $0,$ we have 
 \begin{equation}\tag{2}\label{eqn2}
 \sum_{1 \leq p \leq x}   a_f(p) = \BigO{x^{\frac{155}{156} +\epsilon}}, 
  \end{equation}
 where the sum runs over all prime numbers $p$ $\leq$ $x.$
 \end{prop}

We need a modified version of Proposition 2.2. in the context of almost primes. To this end, we introduce standard functions and apply
Vaughan's identity to obtain this modified version. For any $r\in\N$, let $\Lambda_r$ denote the generalized von Mangoldt function of order $r$. It is defined as
$$\Lambda_r(n) = \sum_{d|n} \mu(d)\left( \mathrm{log} \frac{n}{d}\right)^r.$$
By M\"{o}bius inversion, $$\left( \mathrm{log}{n} \right)^r = \sum_{d|n} \Lambda_r(d).$$ 
For any $r \in \N$, we know that $$\zeta^{(r)}(s)=(-1)^{r-1} \sum_{n=1}^{\infty} (\mathrm{log} n )^r n^{-s}.$$ 
This implies that
$$ \frac{\zeta^{(r)}(s)}{\zeta(s)} = (-1)^r\sum_{n=1}^{\infty} \Lambda_r(n) n^{-s}.$$
The important feature of $\Lambda_r(n)$ is that it vanishes whenever
$n$ has more than $r$ prime factors and so it is useful to detect
when the condition $\omega(n)\leq r$ holds.
We next recall a combinatorial partition of $\Lambda_1(n)$ which is a result of Vaughan ~\cite{Va75} [cf. equation (61) of ~\cite{DI90}].
\begin{prop}
\label{key-identity3}
$$\Lambda_1(n) = \sum_{d|n, d\leq R} \mu(d) \mathrm{log} \frac{n}{d}   
 - \sum_{lm|n, m \leq R, l \leq Q}  \mu(m)\Lambda_1(l).$$  
\end{prop}
We need to generalize this identity to the case $\Lambda_r(n).$ We do this below.

\section{Sums over almost primes and the proof of Proposition \ref{lem1}} 

In Proposition \ref{lem11}, Duke-Iwaniec obtained the estimation for sums over
primes. We modify their proof to show more generally that for any $r$ $\geq$ $1,$
we have Proposition \ref{lem1}.

The following combinatorial partition of $\Lambda_r(n)$ can also be deduced from
[Lemma 1,~\cite{MS02}] with $c(n) = (-1)^r \Lambda_r(n),$
$\tilde{b}(n) = \mu(n),$ $b(n) = 1$ and $a(n) = (\mathrm{log}n)^r.$ 
\begin{prop}
\label{key-identity}
\begin{align*}\tag{*}{\label{eqn * }}
\Lambda_r(n) =& \sum_{d|n, d\leq R} \mu(d) \left(\mathrm{log} \left(\frac{n}{d}\right)\right)^r   
- \sum_{lm|n, m \leq R, l \leq Q}  \mu(m)\Lambda_r(l) \\ &+ \sum_{lm|n,  l \leq Q}  \mu(m)\Lambda_r(l)  + \sum_{lm|n, m > R, l>Q}  \mu(m)\Lambda_r(l).  
\end{align*}
\end{prop}
\begin{proof}
We have for any $y,$
$$ \Lambda_r(n)= \sum_{d|n, d\leq y} \mu(d) \left(\mathrm{log} \left(\frac{n}{d}\right)\right)^r + 
\sum_{d|n, d > y} \mu(d) \left(\mathrm{log} \left(\frac{n}{d}\right)\right)^r $$
Now taking the second sum in this last line we have, by M\"{o}bius inversion,
$$\sum_{d|n, d > y} \mu(d) \left(\mathrm{log} \left(\frac{n}{d}\right)\right)^r
=  \sum_{d|n, d > y} \mu(d) \sum_{c |(n/d)} \Lambda_r(c)$$
$$ = \sum_{cd|n, d > y}  \mu(d)\Lambda_r(c) = \sum_{cd|n, d > y, c>z}  \mu(d)\Lambda_r(c) + \sum_{cd|n, d > y, c \leq z}  \mu(d)\Lambda_r(c), \ \mathrm{for} \ \mathrm{any} \ {z > 0}.$$
Here, $z$ is a parameter chosen optimally in our later estimates. It should not be confused with a complex variable.
Again taking the second sum in the final line we have,  
$$\sum_{cd|n, d > y, c \leq z}  \mu(d)\Lambda_r(c)= \sum_{cd|n,  c \leq z}  \mu(d)\Lambda_r(c) - \sum_{cd|n, d \leq  y, c \leq z}  \mu(d)\Lambda_r(c)$$ 
Putting this together, with $y=R,$ $z=Q,$ we get the required  partition of $ \Lambda_r(n).$
\end{proof}
We obtain the following corollary which is a generalized statement of Proposition ~\ref{key-identity3} to any $r\in \N.$
\begin{cor}
\label{key-identity2}
Suppose $n \in \N$ with $ Q<n \leq QR=X$. Then 
$$\Lambda_r(n) = \sum_{d|n, d\leq R} \mu(d) \left(\log \left(\frac{n}{d}\right)\right)^r   - \sum_{lm|n, m \leq R, l \leq Q}  \mu(m)\Lambda_r(l)$$ 
\end{cor}
\begin{proof}
If $n \leq QR,$ note that the last sum in Proposition \ref{key-identity} is zero. 
If $ n > Q,$ then the third sum in (\ref{eqn * })  $\sum_{lm|n,  l \leq Q}  \mu(m)\Lambda_r(l)  =  \sum_{l|n,  l \leq Q} \Lambda_r(l) \sum_{ m | \frac{n}{l}  } \mu(m)$ is zero since $\frac{n}{l} >1$ because $l \leq Q$ and $n = \frac{n}{l} l > Q.$
\end{proof}
Let $\hat{f_n}=n^{\frac{k-1}{2}} a_n(f)$. Take $b_n = \left(1-\frac{n}{X} \right) \hat{\psi}(n)$, where $\hat{\psi}$ is the Gauss
sum of a Dirichlet character $\psi$ to modulus $c \equiv 0 \pmod {4}$ and $X \geq 2$. Now, let us consider the sum
$$ P(X) = \sum_{n \leq X} b_n \hat{f_n} \Lambda_r(n).$$
We follow closely the method of Duke and Iwaniec \cite{DI90}. 
Similar to the proof of Proposition ~\ref{lem11} in \cite{DI90}, we shall split the second sum in Corollary~\ref{key-identity2}
over the dyadic intervals $L < l \leq 2L$, $M <m \leq 2M$ with $2L \leq Q$ and $2M \leq R$, and we write accordingly
$$\Lambda_r(n) = \Lambda^{*}_R(n) - \sum_L \sum_M \Lambda^{*}_{LM}(n),$$
where $$\Lambda^{*}_R(n)= \sum_{d|n, d\leq R} \mu(d) \left(\mathrm{log} \left(\frac{n}{d}\right)\right)^r,
\Lambda^{*}_{LM}(n)=\sum_{lm|n  , L < l \leq 2L,  M < m \leq 2M} \mu(m) \Lambda_r(l). $$
\begin{lem}
$$P(X) = P_R(X) - \sum_L \sum_M P_{LM}(X) + \BigO{Q^{\frac{k+1}{2}} X^{r\epsilon}},$$
where $$P_R(X) = \sum_{n \leq X} b_n \hat{f}_n \Lambda^{*}_R(n),\quad P_{LM}(X) = \sum_{n \leq X} b_n \hat{f}_n \Lambda^{*}_{LM}(n). $$
\end{lem}
\begin{proof}
The contribution to the sum from $n \leq Q$ is clearly $\BigO{Q^{\frac{k+1}{2}}X^{\epsilon}}$ by a simple application
of Cauchy's inequality as in \cite{DI90}.
\end{proof}
\section{Proof of Proposition \ref{lem1}}
\begin{proof}
We follow \cite{DI90}, but give more details as the proof in \cite{DI90} is terse.
To treat $P_R(X),$ we apply partial summation and (58) of \cite{DI90} which states
$$ \sum_{n \leq X, d | n} \hat{\psi}(n) \hat{f}_n \ll   X^{\frac{k}{2} }\mathrm{log} X ,$$
where the implied constant is independent of $d$.
We deduce $$ P_R(X) \ll R X^{\frac{k}{2}} (\mathrm{log} X)^{r+1}.$$
To treat $P_{LM}(X),$ we follow again \cite{DI90} and split 
$$P_{LM}(X) = P^{\prime}_{LM}(X) + P^{\prime \prime}_{LM} (X)$$
where in the first term, $n = lm$ is squarefree and in the second term $n = lm$ is not
squarefree. Proceeding as in \cite{DI90}, we have (using (7) of \cite{DI90}),
$$P^{\prime}_{LM}(X) \ll LM X^{  \frac{k}{2} - \frac{1}{4} + \epsilon }.$$

For $P^{\prime \prime}_{LM} (X),$ we have (upon using (6) of \cite{DI90}), 

$$P^{\prime \prime}_{LM} (X) \ll \sum_{l,m}  |\mu(m)| |\Lambda_r(l)| (lm)^{\frac{1}{2}} X^{\frac{k-1}{2} +\epsilon},$$
where the sum runs over $l$ and $m$ such that $L \leq l \leq 2L,$ $M \leq m \leq 2M$ and $\mu(lm) = 0.$
Since $\mu(lm)=0 0$, and $l$ is $P_r$, we see that
$(l,m)\neq 1$ in the sum.  We write
$$L < l=p_1 ...p_r < 2L $$
and suppose $p_1|m$.  Thus,
$$P_{LM}^{''} (X) \ll X^{(k-1)/2 + \epsilon}\sum_{L < p_1 < 2L} \sum_{L/p_1 < p_2...p_r < 2L/p_1}
{M^{3/2} \over p_1} \ll LM^{3/2} X^{(k-1)/2+\epsilon}. $$
Hence $$P^{\prime \prime}_{LM} (X) \ll LM^{\frac{3}{2}} X^{\frac{k-1}{2} +\epsilon}.$$
Thus we get the first bound, $$P_{LM} (X) \ll LM X^{  \frac{k}{2} - \frac{1}{4} + \epsilon } + LM^{\frac{3}{2}} X^{\frac{k-1}{2} +\epsilon }. $$  
For the second bound, one appeals to the estimate for the bilinear form:

$$| \sum_{mn \leq X, M < m \leq 2M} a_m b_n \hat{f}_{mn}| \ll ( \sum_{mn \leq 2X} |a_m b_n|^2 )^{\frac{1}{2}} (X^{\frac{1}{2}} M^{-\frac{1}{2}} + X^{\frac{1}{4} } M^{\frac{3}{4}} ) X^{\frac{k-1}{2} +\epsilon}$$
which is (57) in \cite{DI90}. We put $a_m = \mu(m)$ to get
$$| \sum_{mn \leq X, M < m \leq 2M} \mu(m) b_n \hat{f}_{mn}| \ll (M^{-\frac{1}{2}} + M^{\frac{3}{4}} X^{-\frac{1}{4}}) X^{\frac{k+1}{2} +\epsilon}.$$
Combining this with our earlier discussion of $P_{LM}(X)$ finally leads to
$$P_{LM}(X) \ll (M^{-\frac{1}{2}} + M^{\frac{3}{4}} X^{-\frac{1}{4}}) X^{\frac{k+1}{2} +\epsilon}$$
which is valid for any $M, Q, R$ satisfying $1 \leq M \leq R = \frac{X}{Q}.$
Choosing $M = X^{\frac{1}{26}}, Q = X^{\frac{9}{13}}, R = X^{\frac{4}{13}}$ gives
$$P(X) \ll X^{\frac{k+1}{2} - \frac{1}{52} + \epsilon}.$$
The smoothing factor $(1 - \frac{n}{X} )$ is removed as in \cite{DI90}.
Therefore, we obtain that 
$$ \sum_{n \leq X} \hat{\psi}(n) \hat{f}_n \Lambda_r(n) \leq X^{\frac{k+1}{2}-\frac{1}{156}+\epsilon}. $$
Since $\Lambda_r(n)$ vanishes if $n$ has more than $r$ prime factors, we have
by partial summation
$$ \sum_{n \leq X, \omega(n)\leq r} \hat{\psi}(n) a_f(n)  \leq X^{\frac{155}{156}+\epsilon}. $$
Hence $$\sum_{1 \leq n \leq x, \omega(n) \leq r} \hat{\psi}(n) a_f(n) = \BigO{x^{\frac{155}{156} +\epsilon}}.$$
As in \cite{DI90}, we may take $\psi$ to be the principal character
(mod 4), from which the proposition follows.

\end{proof}
\section{Proof of Proposition \ref{mainthm2}}
\begin{proof}
Let $r$ be as in the theorem of \cite{HL97} and $F$ a non-negative smooth
function compactly supported in $(0,1)$ with positive mean value.  
The argument on page 439 of \cite{HL97} shows that 
there is a positive constant $c$ such that
$$\sum_{n \leq Y, \ \omega(n) \leq r } a_f^2(n) F(\frac{n}{Y})\geq 
{cY \over \log Y} + O(Y^{14/15}). $$
Since $F$ is bounded, this means there is a positive constant $c_1$ such that
$$\sum_{n \leq Y, \  \omega(n)\leq r}a_f^2(n) \geq {c_1Y \over \log Y} + O(Y^{14/15}). $$
On the other hand, we have assuming the Ramanujan conjecture, 
$$\sum_{n\leq Y^{\delta}} a_f^2(n) \ll Y^{\delta+\epsilon}. $$
Actually, by the techniques of \cite{murty-murty} and \cite{murty-murty2}, the $\epsilon$ in the
exponent can be removed and one does not need to assume Ramanujan here.
 In any case, for $\delta >0$ and
sufficiently small, we have
$$\sum_{Y^\delta < n < Y, \, n=P_r} a_f(n)^2 \gg {Y \over \log Y}. $$
This completes the proof of Proposition \ref{mainthm2}.
\end{proof}
\section{Proof of Theorem \ref{mainthm1}}
Theorem \ref{mainthm1} follows from the Conjecture \ref{conj}, Proposition \ref{lem1} and Proposition \ref{mainthm2}
but for the convenience of the reader we include a
proof.
\begin{proof}
We choose $\delta$ sufficiently small and show that $a_f(n)$ changes sign for
$n = P_r$ and $x^\delta < n < x$.  Suppose not.  
Without loss of generality, we can assume that $a_f(n)$ are positive for all $n$ in 
the set  $T=\{n : x^\delta < n  \leq x , \omega(n) \leq r\}.$
From Proposition \ref{lem1}, for sufficiently large $x,$ and
sufficiently small $\delta >0$, we have
\\
\begin{equation}\tag{6}\label{eqn6}
\begin{gathered}
\sum_{n \in T} a_f(n)  = \BigO{ x^{\frac{155}{156}+\epsilon } }
\end{gathered}
\end{equation}
Using Conjecture \ref{conj}, 
\begin{equation}\tag{7}\label{eqn7}
\begin{gathered}
\sum_{n \in T} a^2_f(n) = \BigO{ x^{\frac{155}{156}+ \epsilon_0 +\epsilon } }
\end{gathered}
\end{equation}
\\
Replacing $Y$ by $x$ in Proposition \ref{mainthm2}, we get
\\
\begin{equation}\tag{8}\label{eqn8}
\begin{gathered}
\sum_{{\substack{ n = P_r  \\ 1 \leq n \leq x }}} a_f^2(n)  \gg \frac{x}{\log x}.  
\end{gathered}
\end{equation}
Hence for sufficiently large $x$,  we have
\\
\begin{equation}\tag{9}\label{eqn9}
\begin{gathered}
\sum_{n \in T}  a^2_f(n) \gg \frac{x}{\log x}.
\end{gathered}
\end{equation}
\\
We have a contradiction from (\ref{eqn7}) and  (\ref{eqn9}). Thus there is atleast one sign change
of $a_f(n)$  with $n = P_r$  in $(x^\delta, x).$
Thus, there is a sign change in each of the intervals of the form
$(x^{\delta^t}, x^{\delta^{t-1}})$.  The number of such disjoint intervals 
covering $(1,x)$ is clearly $\gg \log x$.   
This completes the proof of our theorem.
\par
The value of $r$ is determined by the results of \cite{HL97}.  In their
paper, the authors suggest that 
by using metaplectic techniques, one can show that $r=4$ 
is permissible. Again based on comments of that paper,
it is possible to sharpen this to $r=3$ by using weighted sieve techniques [cf. see section 3, ~\cite{HL97} for the details].
This seems to be the limit of present-day knowledge.
\end{proof}
\section{Concluding Remarks}
In this section, by assuming a Siegel-type conjecture, we deduce that the
sequence $a_f(p)$, where $p$ varies over primes,
change signs infinitely often. 
\begin{conj}\label{conj2}
If $L(\frac{1}{2}, g, \chi_p) \neq 0,$ then
$|L(\frac{1}{2}, g, \chi_p)| \gg p^{-\epsilon}$ for any $\epsilon$  $>$ $0.$
\end{conj}
\begin{thm}\label{fin-rem}Let $f$ be as in Theorem~\ref{mainthm1}.
Assume that the conjecture \ref{conj2} holds. Then the
sequence that $a_f(p)$, where $p$ varies over primes,
change signs infinitely often. 
\end{thm}
\begin{proof}
Without loss of generality, for sufficiently large $x,$ we can assume that $a_f(p)$ are positive for all $p$ in 
the set $T^{\prime}=\{n : x_0 < p  \leq x \}$  for some natural number $x_0.$ 
From Proposition \ref{lem11}, we have
\\
\begin{equation}\tag{10}\label{eqn10}
\begin{gathered}
\sum_{p \in T^{\prime}} a_f(p) = \BigO{ x^{\frac{155}{156}+\epsilon}}
\end{gathered}
\end{equation}
Using Conjecture \ref{conj2}, and Waldspurger's theorem
\\
\begin{equation}\tag{11}\label{eqn11}
\begin{gathered}
\sum_{p \in T^{\prime} } a_f^2(p)   \gg  x^{1-\epsilon^{\prime}}.  
\end{gathered}
\end{equation}
\\
We have a contradiction from (\ref{eqn10}) and  (\ref{eqn11}). Thus there is atleast one sign change
of $a_f(p)$ for $p$ prime in $(x_0, x].$  As before, this argument 
can be fine-tuned to yield $\gg \log x$
sign changes for $p\leq x$.   
This completes the proof of our theorem.

\end{proof}
These results certainly extend to higher level since both the Duke-Iwaniec
theorem and the Hoffstein-Luo theorem do and our argument goes through.
Finally, we remark that the assumption of 
Ramanujan's conjecture in our main theorem can be relaxed 
somewhat.  
A weaker assumption, namely
$a_f(n) =  \BigO{n^{\alpha}}$ for any $\alpha$ such that $0 < \alpha < \frac{1}{156}$ is sufficient to prove the results.

\subsection*{\sl Acknowledgements}
The first author is grateful to Dr. Narasimha Kumar for helpful discussions.
She would like to thank MPIM (Bonn) and IIT (Hyderabad) for their hospitality.
Part of this work was done when she was a Postdoctoral Fellow at IMSc, Chennai. 
The second author thanks IMSc, Chennai for its kind hospitality and support
while this work was being done.  We thank the referee for very helpful remarks that
considerably enhanced the exposition of this paper.  We also thank
Akshaa Vatwani for spotting several typos.


\begin{thebibliography}{100}
\bibitem{CM06}
              A. C. Cojocaru and  M. R. Murty, 
              An introduction to sieve methods and their applications, 
              London Mathematical Society Student Texts, 66. Cambridge University Press, Cambridge, 2006. 

\bibitem{blomer}  V. Blomer and G. Harcos, Hybrid bounds for twisted $L$-functions, {\sl J. reine angew. Math., \bf 621} (2008), 53-79.



\bibitem{DI90}  W. Duke, William and H. Iwaniec, 
                Bilinear forms in the Fourier coefficients of half-integral weight cusp forms and sums over primes,
                Math. Ann. 286 {\bf 4} (1990), 783--802.            
\bibitem{HL97} J. Hoffstein and W. Luo,  
               Nonvanishing of $L$-series and the combinatorial sieve, 
               With an appendix by David E. Rohrlich. Math. Res. Lett. {\bf4} (1997), 435--444.               
\bibitem{HKKL12} T.A. Hulse, E.M. Kiral, C.I. Kuan and L. Lim,
                 The sign of Fourier coefficients of half-integral weight cusp forms, 
                 Int. J. Number Theory {\bf8} (2012), 749--762. 
\bibitem{HR74}   H. Halberstam and H.-E. Richert, 
                 Sieve methods,
                 Academic Press, London, 1974.
\bibitem{Iwa87} H. Iwaniec, 
                Fourier coefficients of modular forms of half-integral weight, 
                Invent. Math. {\bf87} (1987), no. 2, 385--401.    
\bibitem{Koh10}  W. Kohnen, 
                 A short note on Fourier coefficients of half-integral weight modular forms,
                 Int. J. Number Theory {\bf6} (2010) 1255--1259
\bibitem{Koh92}  W. Kohnen, 
                 On Hecke eigenforms of half-integral weight,
                 Math. Ann. {\bf293} (1992), no. 3, 427--431.
          
\bibitem{kz}  W. Kohnen and D. Zagier, Values of $L$-series of modular forms at the center of the
critical strip, {\sl Inventiones Math., \bf 64} (1981), 175-198.

\bibitem{MM13}  J. Meher and M. R. Murty, 
                Sign changes of Fourier coefficients of half-integral weight cusp forms,
                Int. J. Number Theory {\bf10} (2014), no. 4, 905--914.
                 
\bibitem{Mor73} C.J. Moreno, 
                The Hoheisel phenomenon for generalized Dirichlet series, 
                Proc. Amer. Math. Soc. {\bf40} (1973), 47--51. 
                
                
\bibitem{Mur83} M. R. Murty,
                Oscillations  of Fourier  Coefficients  of Modular  Forms, 
                Math.  Ann. (1983), no. 262, 431--446. 


\bibitem{murty-murty} M. Ram Murty and V. Kumar Murty, Mean values of
derivatives of modular $L$-series, Annals of Math., {\bf 133} (1991), 447-475.

\bibitem{murty-murty2}  M. Ram Murty and V. Kumar Murty, Non-vanishing of $L$-functions and applications, Progress in Mathematics, Volume 157, Birkh\"auser, 1997.
                
                
\bibitem{MS02}  M. R. Murty and A. Sankaranarayanan, 
                Averages of exponential twists of the Liouville function,
                Forum Math {\bf14} (2002), 273--291.
                
\bibitem{Shi73} G. Shimura,
                On modular forms of half integral weight.
                Ann. of Math. (2) {\bf97} (1973), 440--481.  
\bibitem{Va75}  R.C. Vaughan, 
                Mean value theorems in prime number theory, 
                J. Lond. Math. Soc. (2) {\bf10} (1975), 153--162.              

\bibitem{Wa81}  J.-L. Waldspurger,
                Sur les coefficients de Fourier des formes modulaires de poids demi-entier,
                J. Math. Pures Appl. (9) {\bf60} (1981), no. 4, 375--484. 
 
\end{thebibliography}
\end{document}